\documentclass[A4,12pt]{article}

\title{Biochemical pathways simulation}
\author{A. Cangiani}

\newcommand{\nwt}{\newtheorem}
\nwt{theo}{Theorem}
\nwt{lem}{Lemma}
\nwt{prop}{Proposition}
\nwt{defn}{Definition}
\newcommand{\nwc}{\newcommand}
\nwc{\noi}{\noindent}
\nwc{\non}{\nonumber}
\nwc{\p}{\partial}
\nwc{\half}{\frac{1}{2}}
\nwc{\ca}[1]{{\cal #1}}

\nwc{\T}{\textrm}
\nwc{\wh}{\widehat}
\nwc{\wt}{\widetilde}
\nwc{\bfm}[1]{\mbox{\boldmath $#1$}}
\nwc{\vmm}[1]{\vspace{#1 mm}}
\nwc{\hmm}[1]{\hspace{#1 mm}}
\nwc{\Natural}{\mathbb{N}}

\usepackage{amsmath,amssymb}
\usepackage{bm,bbm}

\usepackage{graphicx,graphics,rotating}
\usepackage{overpic}
\graphicspath{{./figures/}}
\usepackage{chemarr}

\usepackage{citesort}

\usepackage{multicol,multirow}

\usepackage{comment}

\usepackage{easyeqn}
\usepackage{easyvector}

\newcommand{\Xvar}{X}                   %% variables
\newcommand{\xvar}{x}                   %% concentration variables
\newcommand{\MATSYMBOL}[1]{\mathsf{#1}}
\newcommand{\Lmat}{\MATSYMBOL{L}}
\newcommand{\Rmat}{\MATSYMBOL{R}}
\newcommand{\Dmat}{\MATSYMBOL{D}}
\newcommand{\Dkmat}{\MATSYMBOL{Dk}}
\newcommand{\Xmat}{\MATSYMBOL{X}}
\newcommand{\kv}{\mathbf{k}}          %% vector of kinetic constants 
\newcommand{\xv}{\mathbf{x}}          %% vector of concentraions
\newcommand{\xxv}{\mathbf{xx}}          %% vector of concentraions
\newcommand{\Xv}{\mathbf{X}}          %% vector of concentraions
\newcommand{\Rv}{\mathbf{R}}          %% vector of concentraions
\newcommand{\fn}{f_n}                 %% single ODE functio
\newcommand{\fv}{\mathbf{f}}          %% ODE function
\newcommand{\lv}{\mathbf{l}}          %% ODE function

\begin{document}

\maketitle

\begin{abstract}
In this short note we review deterministic simulation of biochemical pathways,
i.e. networks of biochemical reactions obeying the law of mass action.
It is meant as a basis for the MATLAB code, written by the author, which permits
easy input and simulation of general biochemical networks. 
This work was carried out for the  European Project `CardioWorkBench'.
\end{abstract}

\section{Introduction}

These notes give a short guidance to deterministic simulation of biochemical pathways.
We follow the approach introduced by Ullah {\em et al.}~\cite{Ullah06}, basis of
the system biology sbtoolbox~\cite{Schmidt06} freely available 
from~\cite{Schmidt}. 
For a comprehensive list of system biology software see~\cite{SBML}. 
% SBML format    export interface

The aim is to easily input and simulate a general biochemical network of reactions obeying
the law of mass action. The pathway specifications (rates, reactants, and products of each
reaction) have to be easily and flexibly importable. And the evaluation of the resulting  system 
of ODEs expressing  the model's dynamics have to be computationally efficient.  
We optimise previous implementations and show how to efficiently calculate  the ODEs jacobian.
We also mention how to easily include parameters like phases and volume
scalings.  
A user-friendly MATLAB code implementing the algorithm discussed in these notes can be 
downloaded from~\cite{CangiODE}.
%Although not discussed here, the interface for import of SBML models
%The compatible for importing models data from SMBL-enabled modelling packages

\section{Modeling}

We consider simulating a biochemical reaction pathway (network)
assuming that it can be decomposed into unidirectional
elementary reactions and that the law of mass action can be applied
to each reaction.

Let $\Xv=(\Xvar_n)_{N\times 1}$ represent the molecular species, 
$\Rv=(R_m)_{1\times M}$, the elementary reactions, 
$\kv=(k_m)_{1\times M}$  the {\em rate} coefficients.
The set of elementary reactions forming the pathway can be written as
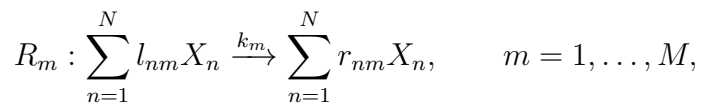
\begin{EQ}\label{eq:elementary_reactions}
R_m:\sum_{n=1}^N l_{nm} \Xvar_n\xrightarrow{k_m}
\sum_{n=1}^N r_{nm} \Xvar_n,\qquad m=1,\dots,M,
\end{EQ}
where the reactions' {\em stoichiometric}
coefficients $l_{nm}$ and $r_{nm}$ are non-negative integers.
%Notice that, since the reactions are elementary,
%one of $l_{nm}$ and $r_{nm}$ must be zero.
The step change in the number of molecules $\Xvar_n$ due to
reaction $R_m$ is given by 
%thus given by
\[
d_{nm}=r_{nm}-l_{nm}.
\]
We collect these coefficients into the input and output
stoichiometric matrices $\Lmat=( l_{nm} )_{N\times M}$ and 
$\Rmat=( r_{nm} )_{N\times M}$, and the step change matrix
$\Dmat=( d_{nm} )_{N\times M}$.

In practice, it is convenient to input the information contained in
the matrices $\Lmat$ and $\Rmat$  by storing  their (few)
non-zero entries values and indices, as it is done in any matrix 
{\em sparse} representation.  This process is described below with an example.
%Let $N_l$ and $N_r$ be, respectively,  the maximum number of reactants
%and products of a single elementary reaction.
%We define two new matrices L^
%%From the pratical point of view,
%%The process of simulating a bichemical pathway can be decomposed into
%%wo tasks: the first is to set up the ODE system by obtaining the matrices
%%$\Lmat$ and $\Rmat$; the second is the actual solution of the system
%%with the given initial values.
%%The letter task will be discussed later on. As for the former,
%%
%%\eqref{eq:odefun}

The law of mass action implies that
the dynamics of the molar concentrations $\xvar_n=[\Xvar_n]$,
$n=1,\dots,N$ is  described by a  system of ODEs, each called a {\em rate equation},
that we can easily
write down in terms of the reactions' coefficients.
The set of molar concentrations 
$\xv=(\xvar_n)_{N\times 1}$ must satisfy the system of first order
ODEs
%system describing the dynamics of each concentration as
\begin{EQA}[l]
\label{eq:odesys}
\dot{\xv}=\fv(\xv)\quad\T{with}\quad \fv=(\fn)_{N\times 1}\quad\T{and}\\
\label{eq:odefun}
\fn(\xv)=\sum_{m=1}^M d_{nm} \left(k_m \prod_{i=1}^N x_i^{l_{im}}\right).
\end{EQA} 
The parameter $d_{nm}$ is different from zero if the $n$-th variable is
involved in the $m$-th reaction, the rate of which is given by the expression 
in parentheses.
Note that all this can be generalized to include more complex dynamics
by allowing the stoichiometric coefficients to be non-negative real 
numbers (Generalized Mass Action), see~\cite{Voit}.

\

{\bf Example.} Consider the basic enzyme-kinetic reaction
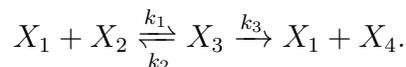
\begin{EQ}
\nonumber
\Xvar_1+\Xvar_2 \xrightleftharpoons[k_2]{k_1}
\Xvar_3\xrightarrow{k_3}
\Xvar_1+\Xvar_4.
\end{EQ}
%where $\Xvar_1,\Xvar_2,\Xvar_2$, and $\Xvar_4$ are the number of
%enzyme, substrate, complex, and product, respectively.
We decompose the reaction pathway into
the elementary reactions
\begin{EQ}[l]
\nonumber
R_1:\Xvar_1+\Xvar_2\xrightarrow{k_1} \Xvar_3
\\
\nonumber
R_2:\Xvar_3 \xrightarrow{k_2} \Xvar_1+\Xvar_2
\\
\nonumber
R_3:\Xvar_3 \xrightarrow{k_3} \Xvar_1+\Xvar_4
\end{EQ} 
To this sequence of reactions we associate the following coefficient matrices
used to compute~\eqref{eq:odefun}:
\begin{EQ}
\nonumber
\Dmat=\left[
\begin{array}{ccc}
-1 & 1  & 1 \\
-1 & 1  & 0 \\
 1 & -1 & -1 \\
 0 & 0  & 1 
\end{array}
\right],\quad
\Lmat=\left[
\begin{array}{ccc}
1 & 0 & 0 \\
1 & 0 & 0 \\
0 & 1 & 1 \\
0 & 0 & 0 
\end{array}
\right].
\end{EQ}
These matrices can be  obtained from a compact representation,
that facilitates the input of the model data.
For instance, the matrix $\Lmat$ specifying the reactants, can be 
computed from 
\begin{EQ}
\Lmat^\T{i}=
\left[
\begin{array}{ccc}
1 & 3 &  3 \\
2 & 0 & 0 
\end{array}
\right],\qquad
\Lmat^\T{v}=
\left[
\begin{array}{ccccc}
1 & 1 &  1 \\
1 & 0 & 0
\end{array}
\right].
\end{EQ}
The first row of $\Lmat^\T{i}$ (the superscript ``$\T{i}$'' stands for $\T{i}${\em ndex}) 
tells us that the first reaction involves the first
and second reactant, and so on. The matrix $\Lmat^\T{v}$ 
(``$\T{v}$'' for  $\T{v}${\em alue}) collects 
the related non-zero stoichiometric coefficients $l_{nm}$.
All remaining entries are filled with zeros.
Similarly, we decompose the matrix $\Rmat$ which contains the information 
about the reactions' products into two matrices $\Rmat^i$ and $\Rmat^v$.

In realistic biochemical pathways, there are many reactants and reactions
involved but few reactants acting on the single elementary reaction.
Consequently, most of the entries of the matrices 
$\Lmat$, $\Rmat$ (and $\Dmat$) are zero. The compact (sparse) 
representation above represents a good compromise between user friendliness 
and computational efficiency.

\section{Implementation and extensions}

The process of simulating a biochemical pathway can be decomposed into
two tasks: the first is to set up the ODE system by obtaining the matrices
of stoichiometric coefficients $\Lmat$ and $\Rmat$, and the vector of  reaction rates $\kv$; 
the second is the actual numerical solution (time-stepping) of the
system with the given concentrations initial values.
The computational cost of time-stepping is dominated by the evaluation of the ODEs' function, thus it is crucial that this is done efficiently.

The implementation of $\fv(\xv)$ presented in~\cite{Ullah06} uses the matrices
$\Lmat$ and $\Dmat$. Letting $\xxv=(\xv,\dots,\xv)_{N\times M}$ be a matrix
whose $M$-columns are 
copies\hspace{-1 mm}~\footnote{
in MATLAB, this is achieved by the command xx=repmat(x,1,M)
%%%%%%%%%%%%%%%%%%%%%%%
}
of $\xv$, the expression in~\eqref{eq:odefun} can be evaluated as follows:
\begin{EQ}
%\fv(\xv)=\Dmat(\kv\bullet\:\T{prod}(\xxv\bullet \Lmat,1))^T,
\fv(\xv)=\Dmat\ast(\kv\,.\hspace{-1 mm}
\ast\T{prod}(\xxv\, .^\wedge \Lmat,1))',
\end{EQ}
where $\ast$ denotes matrix multiplication, $.\hspace{-.3 mm}\ast$
and $.^\wedge$  entry-wise matrix multiplication and exponentiation,
$\T{prod}(\cdot,1)$ multiplication along columns, and $'$ transposition
(these notational conventions are in accordance with MATLAB's syntax).

In~\cite{Ullah06}, it is also proposed a code  that  speeds-up the evaluation
of $\fv(\xv)$ by avoiding multiplications by zero entries.
This can be achieved more simply by using the compact representation above,
and, if sparse matrices representation is implemented, by creating $\Dmat$
as a sparse matrix. A further straightforward computational optimization is obtained  
by multiplying 
each column of $\Dmat$ with the corresponding rate constant in $\kv$.
Denoting the result of such multiplication by $\Dkmat$, the
computation of $\fv(\xv)$ is reduced to the following steps:
\begin{itemize}
\item[1.]
using $\Lmat^\T{i}$, form the matrix
$\qquad
\Xmat=
\left[
\begin{array}{ccc}
x_1 & x_3 &  x_3 \\
x_2 & 0 & 0 
\end{array}
\right];
$
\item[2.]
evaluate $\fv(\xv)$ as:
$\qquad
\fv(\xv)=\Dkmat \ast
\T{prod}(\Xmat\, .^\wedge \Lmat^\T{v},1)',
$
\end{itemize}
with the convention that the empty product $0^0=1$.

Biochemical models often include parameters like, for instance, phases and volume 
scalings. 
It is useful to implement such parameters separately from the
stoichiometric coefficients. This is easily done  by multiplying the appropriate entry
of $\Dkmat$ by the given parameter.

As the ODE system generated by biochemical pathways is non-linear, it is also useful
to have a routine that evaluates the jacobian $J_\fv(\xv)$ of the ODE system
function $\fv$.
The entries of the jacobian are given by:
\begin{EQ}\label{eq:jacob}
\frac{\partial f_n}{\partial x_j}(\xv)=\sum_{m=1}^M dk_{nm}
%l_{jm} x_j^{\max(l_{jm}-1,0)}\prod_{\stackrel{i=1}{\T{\tiny{$i\neq j$}}}}^N x_i^{l_{im}}.
\left(l_{jm} x_j^{\max(l_{jm}-1,0)}\prod_{i\neq j} x_i^{l_{im}}\right),
\qquad j,n=1,\dots,N.
\end{EQ}

We proceed with the evaluation of the jacobian processing it by columns, i.e. 
evaluating $\frac{\partial \fv}{\partial x_j}$ for $j=1,\dots,N$.

Let $j\in \{1,\dots,N\}$ be given. The evaluation of~\eqref{eq:jacob}
is simplified  by limiting the summation to the reactions involving $\xvar_j$, and
by limiting the multiplications to the reactants involved in the given reaction.
This latter simplification is already embedded  in our compact representation.
As for the summation, we define reduced 
matrices $\Dkmat_{\xvar_j}$, $\Xmat_{\xvar_j}$, and $\Lmat_{\xvar_j}^\T{v}$ obtained from
the corresponding matrices $\Dkmat$, $\Xmat$, and $\Lmat^\T{v}$ by considering only those columns associated to reactions that have $\xvar_j$ among the input reactants.

We need to derive the entries of $\Xmat_{\xvar_j}.^{\wedge}\Lmat^\T{v}$ with respect to the variable $\xvar_j$.
To this end, we define a new row-vector $\lv_{\xvar_j}$ collecting the values $l_{jm}$ 
with $m$ such that $\xvar_j$ is a reactant of reaction $m$, and a new matrix
$\widetilde{\Lmat}_{\xvar_j}^\T{v}$ obtained from $\Lmat_{\xvar_j}^\T{v}$ by replacing 
the entries indexed $jm$ with  $l_{jm}-1$.
Notice that all such matrix manipulations are made particularly easy by MATLAB built-in 
vector/matrix manipulation  tools.
In this way, we can express (and calculate) the $j$-th column of the jacobian matrix in  compact  form as:
\begin{EQ}
\frac{\partial \fv}{\partial x_j}(\xv)=
\Dkmat_{\xvar_j}\ast
\left(
\lv_{\xvar_j}\,.\hspace{-1 mm}\ast  
\T{prod}(\Xmat_{\xvar_j}\, .^\wedge \widetilde{\Lmat}_{\xvar_j}^\T{v})
\right)'.
\end{EQ}

Finally, let us discuss the problem of numerically solving~\eqref{eq:odesys}.
The main characteristics of such system of ODEs are the following.
Unless all elementary reactions are of $0$-th or $1$-st order, the system is 
nonlinear. Moreover, the system is  generally {\em stiff}. 
Thus, it is compulsory to consider variable time-stepping and employ stiff ODE solvers.
Standard ODE solvers packages include robust stiff solvers (for instance, MATLAB's
ode15s), that are fast enough if only a few simulations are needed.
If, on the other hand, a large number of simulations is required as
in parameter estimation, then it may be preferable to code {\em ad hoc} solvers
that take into account the peculiarities of the ODE systems generated  by biochemical pathways.

\nocite{*}
\bibliographystyle{abbrv}
\addcontentsline{toc}{chapter}{Bibliography}
%\bibliography{biosym}

\end{document}